\documentclass[11pt]{article}

\usepackage{amssymb,amsmath,graphicx,amsthm}
\usepackage{sectsty}

\usepackage[english,francais]{babel}

\newtheorem{theorem}{Theorem}[section]
\newtheorem{lemma}[theorem]{Lemma}
\newtheorem{e-proposition}[theorem]{Proposition}

\newtheorem{e-definition}[theorem]{Definition\rm}

\def\og{\leavevmode\raise.3ex\hbox{$\scriptscriptstyle\langle\!\langle$~}}
\def\fg{\leavevmode\raise.3ex\hbox{~$\!\scriptscriptstyle\,\rangle\!\rangle$}}

\def\ie{\emph{i.e. }}

\def\cf{\emph{cf. }}

\def\eps{\epsilon}
\def\rr{\mathbb{R}}

\def\id#1{\mathfrak{#1}}
\def\jo#1{\mathcal{#1}}

\def\limm#1{\textrm{\raisebox{.5ex}{\mbox{$\underset{#1}{\lim}$}}} \:}
\def\lims#1{\textrm{\raisebox{.5ex}{\mbox{$\underset{#1}{\limsup}$}}} \:}
\def\supp#1{\textrm{\raisebox{.5ex}{\mbox{$\underset{#1}{\sup}$}}} \:}

\def\nr#1{\left\| #1 \right\|}
\def\somme#1#2{\overset{#2}{\underset{#1}{\sum}}}

\def\ssi{\Leftrightarrow}
\def\imp{\Rightarrow}

\def\diam{\mathrm{Diam}\,}
\def\dim{\mathrm{dim}\,}

\def\wdm{\mathrm{wdim} \,}
\def\einj{ \!\! \textrm{ \mbox{ $ {}^{\eps} \!\!\!\!\!\!\! \hookrightarrow$ } } \!\!\! }

\def\meas{\mathrm{amp}}

\begin{document}

\sectionfont{\normalsize}
\selectlanguage{english}
\centerline{ {\large On a H{\"o}lder covariant version of mean dimension} }

\selectlanguage{english}
\centerline{Antoine Gournay}
\centerline{ \footnotesize gournay@math.kyoto-u.ac.jp}

{\center Department of Mathematics, Faculty of Science \\
 Kyoto University,\\
 Kyoto 606-8502, Japan

~}

\medskip

\begin{abstract}
\selectlanguage{english}
Let $\Gamma$ be a infinite countable group which acts naturally on $\ell^p(\Gamma)$. We introduce a modification of mean dimension which is an obstruction for $\ell^p(\Gamma)$ and $\ell^q(\Gamma)$ to be H{\"o}lder conjugates.  

\vskip 0.5\baselineskip

\selectlanguage{francais}
\noindent{\bf R\'esum\'e} \vskip 0.5\baselineskip \noindent
{\bf Sur une modification H{\"o}lder covariante de la dimension moyenne. }
Soit $\Gamma$ un groupe d{\'e}nombrable infini qui agit naturellement sur  $\ell^p(\Gamma)$. Nous introduisons une obstruction, proche de la dimension moyenne, au fait que  $\ell^p(\Gamma)$ et $\ell^q(\Gamma)$ soit H{\"o}lder conjugu{\'e}s.

\end{abstract}

\selectlanguage{english}

\section{Introduction}
\label{}
The Lipschitz and uniform classification of Banach space and their spheres is well-understood (see \cite{BL}), few things are known on their H{\"o}lder-continuous classification. By making hypothesis of a more dynamical nature on the map, we present here some obstructions on the existence of equivariant H{\"o}lder-continuous maps between balls of Banach spaces. 
\par Let $1 \leq p \leq \infty$ and $\Gamma$ be an infinite countable group. Mean dimension (which is only defined when $\Gamma$ is amenable) is a topological invariant of metric spaces with group action introduced in \cite{Gro}. It turns out to be inappropriate to tell apart the dynamical systems given by the unit balls of $\ell^p(\Gamma; \rr)$ when endowed with the natural action of $\Gamma$ and the product topology (induced from $[-1,1]^\Gamma$). Recall that $\ell^p(\Gamma; \rr)$ is the subset of the set of functions $\Gamma \to \rr$ such that $\sum_{\gamma \in \Gamma} |f(\gamma)|^p$ is finite. For $f$ such a function and $\gamma \in \Gamma$ then the natural action is defined by $\gamma f (\cdot) = f(\gamma^{-1}\cdot)$. Note that for the natural metric coming from the $\ell^p$ norm the action of $\Gamma$ is isometric, whereas for any metric compatible with the product topology it is not.
\par When $p<\infty$ these dynamical systems are topologically equivalent, as can be seen by using the Mazur map (see below). We briefly recall that mean dimension itself turns out to be insufficient as all elements are sent to $0$ by the action. By \cite[{\S}5]{LW} the mean dimension is $0$. In \cite[{\S}1.6]{Gro}, an alternative construction of mean dimension is attempted, but to no avail (see \cite{tsu}). The one we present here works as long as metrics of type \eqref{defmet} are used. Introducing a variant of mean dimension will enable us to distinguish these dynamical systems without requiring amenability of $\Gamma$.
\par As we are speaking of H{\"o}lder maps, it is important to specify the metrics used. Consider on $B^{\ell^p}_1 := \{ x \in \ell^p(\Gamma;\rr) | \nr{x}_{\ell^p} \leq 1 \}$, the metric $d^{\{F_i\},\{a_i\}}$ defined as follows. Let $\{F_i\}_{i \geq 1}$ be a sequence of finite sets such that $\cup F_i = \Gamma$. Let $\{a_i\}_{i \geq 1}$ be a summable sequence of positive real numbers normalized so that  if $\sigma_k = \sum_{i \geq k} a_i$, then $\sigma_1 =1$. Let $r\in [1,\infty)$. Recall that $\nr{x_1-x_2}^r_{\ell^r(F_i)} = \sum_{\gamma \in F_i} |x_1(\gamma)-x_2(\gamma)|^r$. Now, let 
\begin{equation}\label{defmet}
d^{\{F_i\},\{a_i\},r}(x_1,x_2) = \Big( \somme{i \geq 1}{} a_i  \frac{1}{|F_i|} \nr{x_1-x_2}^r_{\ell^r(F_i)} \Big)^{1/r}.
\end{equation}
We will allow $r=\infty$ in order to consider the metric $d^{\{F_i\},\{a_i\},\infty}(x_1,x_2) =  \sum a_i \nr{x_1-x_2}_{\ell^\infty(F_i)}$. Usually, the $a_i$ are chosen to be $2^{-i-1}$. The normalization $\sigma_1=1$ is made so that the diameter of the unit ball remains less than $2$. Notice that \mbox{$d^{\{F_i\},\{a_i\}}(x_1,x_2) \leq 2 \sigma_{k+1}$} if $x_1$ and $x_2$ start to differ on $F_{k+1}$, but are identical when restricted to the $F_j$ for $j\leq k$. The topology coming from these metrics is the same as the one induced by the product topology (or the weak-$*$ topology).
\par The result presented here holds if one replaces $\rr$ and $|\cdot|$ by a finite-dimensional $V$ with a norm $\nr{\cdot}$. However, this makes calculations in Lemma \ref{mpq} slightly more tedious, and we will avoid this unnecessary complication. 
\begin{theorem} \label{tprof}
Let $\Gamma$ be an infinite countable group. Let $p,q \in [1,\infty[$. Let $B^{\ell^p}_1$ and $B^{\ell^q}_1$ be unit balls of $\ell^p(\Gamma;\rr)$ and $\ell^q(\Gamma;\rr)$ respectively, both endowed with the action of $\Gamma$ and metrics of the type \eqref{defmet}. When $q > p$, there is no $\Gamma$-equivariant Lipschitz homeomorphism $f: B^{\ell^p}_1 \to B^{\ell^q}_1$. Moreover, if $\alpha q > p $ there is no $\Gamma$-equivariant map $f: B^{\ell^p}_1 \to B^{\ell^q}_1$ which is H{\"o}lder continuous of exponent $\alpha$.
\end{theorem}
\par For certain metrics, the Mazur map is a limit case for the above theorem. Indeed, let $\phi: \ell^p \to \ell^q$ be defined for $x \in \ell^p(\Gamma,\rr)$ by $\phi(x)(\gamma) = \big( x(\gamma) \big)^{p/q}$.
Then, $\phi(B^{\ell^p}_1) = B^{\ell^q}_1$ and $\phi$ is a $\Gamma$-equivariant homeomorphism which is H{\"o}lder continuous of exponent $\alpha = p/q$ if the metric $d^{\{F_i\},\{a_i\},p}$ is used on the source and $d^{\{F_i\},\{a_i\},q}$ on the image (\cf \cite[{\S}9.1]{BL}). Recall (\cf \cite[thm 2.3]{BL}) that any uniformly continuous map (for the norm topology) from $B^{\ell^p}_1 \to \ell^q$ is uniformly approximated by H{\"o}lder continuous maps.
\par The obstruction that is exhibited here is very crude. For example, $\Gamma$ is not required to be amenable and the group structure is almost secondary. In fact, all we really need to know is whether the quantity covariant under H{\"o}lder-continuous equivariant maps is $0$ or $>0$.

\section{Dynamical covariant}
\par We will be using maps that are close to be injective.
\begin{e-definition}
Let $Y$ be a topological space. Let $(X,d)$ be a metric space and  $\eps \in \rr_{>0}$. Then we say a map $f:X \to Y$ is an $\eps$-embedding if it is continuous and the diameter of the fibers is less than $\eps$, \ie $\forall y \in Y, \diam f^{-1}(y) \leq \eps$. 
\end{e-definition}
When one is interested in compression algorithms, it is not only important that the compression map has ``small'' fibers (so that not too much data is lost) but also has an image which is ``small'' in some sense (so that the compression is effective). The next definition attempts to capture the relation between these two features.
\begin{e-definition}
Let $(X,d)$ be a metric space. Call $\wdm_\eps (X,d)$ the smallest integer $k$ such that there exists an $\eps$-embedding $f:X \to K$ where $K$ is a $k$-dimensional polyhedron.
\[
\wdm_\eps (X,d) = \inf_{X \einj K } \dim K.
\]
\end{e-definition}
The function $\eps \mapsto \wdm_\eps (X,d)$ is always decreasing. We can now make the construction that will detect obstructions.
\begin{e-definition}
Let $\id{f} : \jo{P}_{finite}(\Gamma) \to \rr_{> 0}$ a decreasing function (\ie $A \subset B \imp \id{f}(A) \geq \id{f}(B)$), then the asymptotic measure with profile $\id{f}$ will be defined for an increasing sequence of finite sets $\Omega_i$ (whose size tends to $\infty$) by
\[
\meas_\id{f} (X:\{\Omega_i\}) = \lims{i \to \infty} \frac{\wdm_{\id{f}(\Omega_i)} (X,d_{\Omega_i}) }{|\Omega_i|} \quad \in [0,+\infty].
\]
where $d_{\Omega}$ is the dynamical distance: $d_{\Omega}(x,x') = \supp{\gamma \in \Omega} d(\gamma x, \gamma x')$.
\end{e-definition}
The quantity obtained certainly lacks topological invariance and depends on the sequence chosen, but it is nevertheless useful.
\begin{e-proposition} \label{minv}
Let $f:(X,d) \to (X',d')$ be a $\Gamma$-homeomorphism of metric spaces endowed with an action of $\Gamma$ whose modulus of continuity is $\omega_f$. If there exists an increasing function $\psi: \rr_{\geq 0} \to \rr_{\geq 0}$ such that $\limm{\eps \to 0} \psi(\eps) =0 $, and that $\omega_f(\eps) \leq \psi(\eps) $, then
\[
\meas_{\psi \circ \id{f}} (X':\{\Omega_i\}) \leq \meas_{\id{f}} (X:\{\Omega_i\}).
\]
\end{e-proposition}
\begin{proof}
The inequality
\[
d'(f(x),f(y)) \leq \omega_f (d(x,y)) \leq \psi(d(x,y))
\]
passes to $d_\Omega$ and $d'_\Omega$. This yields
\[
\wdm_{\psi(\eps)} (X',d'_{\Omega_i}) \leq \wdm_{\eps} (X,d_{\Omega_i}),
\]
by composing any $\eps$-embedding on $X$ by $f^{-1}$. Choosing $\eps = \id{f}(\Omega_i)$ and passing to the limit gives the claimed inequality.
\end{proof}
Denote by $B^{\ell^p(n)}_r$ the subset of $\rr^n$ of $\ell^p$ norm less than $r$. The definition of $\meas_{\id{f}}$ is made to use the following result from \cite[corollary 1.1]{tsu} or \cite[proposition 1.3]{moi}.
\begin{e-proposition}\label{p1}
For all $\eps \in \rr_{>0}$ and all integers $n\geq 1$ and $0\leq r \leq n-1$, 
\[
\wdm_\eps (B_1^{\ell^p(n)},\ell^\infty) = \left\{
  \begin{array}{llrcl}
    0  & \textrm{if} &             2 \leq & \eps & \\
    r  & \textrm{if} & 2(r+1)^{-1/p} \leq & \eps & < 2r^{-1/p} \\
    n  & \textrm{if} &                    & \eps &< 2n^{-1/p}
  \end{array}. \right.
\]
\end{e-proposition}

\section{Evaluation of $\meas_\id{f}$}
\begin{lemma} \label{mpq}
Let $\Gamma$ be an infinite countable group, let $\{\Omega_i\}$ an increasing sequence of finite subsets, let $p \in [1,\infty[$, let $B^{\ell^p}_1$ be the unit ball of $\ell^p(\Gamma;\rr)$ endowed with a metric of the type \eqref{defmet} and of the natural action of $\Gamma$. For $c \in \rr_{>0}$ and $q \in [1,\infty[$, define the profile  $\id{f}_{c,q}(\Omega) = c |\Omega|^{-1/q}$. Then 
\[
\meas_{\id{f}_{c,q}} (B^{\ell^p}_1:\{\Omega_i\}) \left\{
  \begin{array}{ll}
  > 0           & \textrm{if } p \geq q \\
   = 0          & \textrm{if } p<q
  \end{array}
  \right.
\]
\end{lemma}
\begin{proof}
The arguments that will be used are independent of the choice of the parameters entering in the definition of the metric \eqref{defmet}, which is surprising. We shall abbreviate $d' = d^{\{F_i\},\{a_i\},r}$
\par We will exhibit a lower and an upper bound. Both use proposition \ref{p1} to evaluate $\wdm$ of a ball in $(\rr^n,l^p)$. This proposition can be rewritten as: for $0 \leq k < n$  
\[
k < (\eps/2)^{-p} \leq k+1 \ssi \wdm_\eps (B^{\ell^p(n)}_1,\ell^\infty) =k.
\]
Consequently,
\[
\begin{array}{c}
  \min \bigg( \Big(\frac{\eps}{2}\Big)^{-p}-1,n \bigg) \leq \wdm_\eps (B^{\ell^p(n)}_1, \ell^\infty) < \min \bigg( \Big( \frac{\eps}{2} \Big)^{-p},n \bigg).
\end{array}
\]
This said, we start with the lower bound. This is pertinent only if $p \geq q $. The argument reduces to the construction of an injection which increases distances. Let $\sigma' = \sum_{i\geq 1} a_i/|F_i|$ for $r<\infty$ and when $r=\infty$ take $\sigma'=\sigma_1=1$. Then, for $g \in F_1$, $(B_{\sigma'}^{\ell^p(\Omega_i^{-1} g; \rr)},\ell^\infty)$ injects in $(B^{\ell^p}_1,d'_{\Omega_i})$ by defining the function to be $0$ outside $\Omega_i^{-1} g$. This gives
\[
\begin{array}{c}
\forall k \geq 1, \quad \min \big( ( \eps / 2\sigma' )^{-p}, |\Omega_i| \big) \leq \wdm_\eps (B_{\sigma'}^{\ell^p(\Omega_i^{-1} g, \rr)},\ell^\infty) \leq \wdm_\eps (B^{\ell^p}_1,d'_{\Omega_i}).
\end{array}
\]
Next, dividing both sides by $|\Omega_i|$, taking $\eps = c |\Omega_i|^{-1/q}$ and passing to the limit yields 
\[ 
\lims{i \to \infty} \min \big( ( 2\sigma' / c )^p |\Omega_i|^{\frac{p}{q}-1} , 1 \big) \leq \meas_{\id{f}_{c,q}} (B^{\ell^p}_1,d').
\]
As $p \geq q$ the left-hand side is greater than $\min \big( ( 2\sigma' / c )^p , 1 \big)$. 
\par The upper bound, which is of interest only if $p \leq q$, is obtained by looking at the restriction of $(B^{\ell^p}_1,d'_{\Omega_i})$ to $(B^{\ell^p(\Omega_i^{-1} F_k, \rr)}_1,\ell^\infty)$. If $\eps > \sigma_{k+1}$, this gives a $2\eps$-embedding for $(B^{\ell^p}_1,d'_{\Omega_i})$ given any $\eps$-embedding $f$ of $(B^{\ell^p(\Omega_i^{-1} F_k, \rr)}_1,\ell^\infty)$. Indeed, the $\ell^\infty$ norm is bigger than $d'_{\Omega_i}$ on the restricted set and the restriction map is such that if itself  and $f$ both have fibers of size $\eps$ then the composition has fibers of size at most $2\eps$. Thus, when $\eps > \sigma_{k+1}$, 
\[
         \wdm_{2\eps} (B^{\ell^p}_1,d'_{\Omega_i})   \leq  \min( \eps^{-p}, |\Omega_i^{-1} F_k|) .
\]
For $i$ fixed, take $\eps(i) = c |\Omega_i|^{-1/q}$ and $k(i)$ such that $\eps(i) > \sigma_{k(i)+1} $. In particular, when $i \to \infty $, $\eps(i) \to 0$ and $k(i) \to \infty$. After division by $|\Omega_i| $ and taking the limit, we see that
\[
\meas_\id{f_{c,q}} (B^{\ell^p}_1,d') \leq \lims{i \to \infty} \min \bigg( \big(\tfrac{2}{c} \big)^{p} |\Omega_i|^{p/q-1}, \frac{|\Omega_i^{-1} F_{k(i)}|}{ |\Omega_i|} \bigg).
\]
When $p<q$ this tends to $0$. If $p=q$, things behave differently: the right-hand term is only bounded by $\big(\tfrac{2}{c} \big)^{p}$.
\end{proof}
\section{Proof of the theorem}
If $f: X \to X'$ is H{\"o}lder continuous of exponent $\alpha \in (0,1)$ or Lipschitz (this corresponds to $\alpha =1$), $\exists c' \in \rr_{>0}$ such that the modulus of continuity $\omega_f$ is bounded by $\psi(\eps) = c' \eps^\alpha$. Thus, proposition \ref{minv} gives that $\meas_{\psi \circ \id{f}} B^{\ell^q}_1 \leq \meas_\id{f} B^{\ell^p}_1$, in other words, that $\meas_{\id{f}_{c'c^{\alpha},r/\alpha}} B^{\ell^q}_1 \leq \meas_{\id{f}_{c,r}} B^{\ell^p}_1$. Taking $p=r$, and thanks to lemma \ref{mpq}, we must have that $\alpha q \leq r = p$ to avoid a contradiction.


\begin{thebibliography}{00}

\bibitem{BL}
Y.~Benyamini and J.~Lindenstrauss.
\newblock {\em Geometric nonlinear functional analysis. {V}ol. 1}, volume~48 of
  {\em American Mathematical Society Colloquium Publications}.
\newblock American Mathematical Society, Providence, RI, 2000.

\bibitem{moi}
A.~Gournay.
\newblock Width of $\ell^p$ balls.
\newblock {\em arXiv:0711.3081v2 or http://hal.archives-ouvertes.fr/hal-00189043/fr/}, page~16, 2008.

\bibitem{Gro}
M.~Gromov.
\newblock Topological invariants of dynamical systems and spaces of holomorphic
  maps. {I}.
\newblock {\em Math. Phys. Anal. Geom.}, 2(4):323--415, 1999.

\bibitem{LW}
E.~Lindenstrauss and B.~Weiss.
\newblock Mean topological dimension.
\newblock {\em Israel J. Math.}, 115:1--24, 2000.

\bibitem{tsu}
M.~Tsukamoto.
\newblock Macroscopic dimension of the {$l\sp p$}-ball with respect to the
  {$l\sp q$}-norm.
\newblock {\em J. Math. Kyoto Univ.}, 48(2):445--454, 2008.

\end{thebibliography}
\end{document}